\documentclass[11pt,reqno]{amsart}
\usepackage{diagrams}
\usepackage{amssymb}

\usepackage[T1]{fontenc}
\usepackage[utf8]{inputenc}



\numberwithin{equation}{section}

\theoremstyle{definition}
\newtheorem{thm}{Theorem}[section]
\newtheorem{cor}[thm]{Corollary}

\newtheorem{defi}[thm]{Definition}
\newtheorem{rem}[thm]{Remark}
\newtheorem{note}[thm]{Notation}

\DeclareMathOperator{\msp}{\mathrm{Spec}}

\DeclareMathOperator{\m}{\mathcal{M}}

\DeclareMathOperator{\mo}{\mathcal{O}}

\newcommand{\mr}[1]{\mathrm{#1}}

\newcommand{\mc}[1]{\mathcal{#1}}
\newcommand{\ov}[1]{\overline{#1}}

\begin{document}

\title[$C_1$ conjecture in mixed characteristic]{A pathological case of the $C_1$ conjecture in mixed characteristic}

\author[I. Kaur]{Inder Kaur}

\address{ Instituto Nacional de Matem\'{a}tica Pura e Aplicada, Estr. Dona Castorina, 110 - Jardim Bot\^{a}nico, Rio de Janeiro - RJ, 22460-320,  Brazil
}

\email{inder@impa.br}

\subjclass[2010]{Primary $14$D$20$, $14$J$60$, $14$G$05$, $14$M$22$ Secondary $14$L$24$, $14$D$22$}


\date{\today}

\begin{abstract}
 Let $K$ be a field of characteristic $0$. 
Fix integers $r,d$ coprime with $r\geq 2$.
Let $X_K$ be a smooth, projective, geometrically connected curve of genus $g\geq 2$ defined over $K$. 
Assume there exists a line bundle $\mc{L}_K$ on $X_K$ of degree $d$. 
In this article we prove the existence of a stable locally free sheaf on $X_K$ with rank $r$ and determinant $\mc{L}_K$. 
This trivially proves the $C_1$ conjecture in mixed characteristic for
the moduli space of  stable locally free sheaves of fixed rank and determinant over a smooth, projective curve.
\end{abstract}

\maketitle

\section{Introduction}

A field $L$ is said to be $C_{1}$ if any hypersurface in $\mathbf{P}^{n}_{L}$ of degree $d \leq n$ has a rational point. The Lang-Manin-Koll\'ar conjecture states that a smooth, proper, separably rationally connected variety over a $C_1$ field has a rational point. Let $K$ be the fraction field of a Henselian discrete valuation ring with algebraically closed residue field denoted $k$. By {\cite[Theorem $14$]{L}}, $K$ is a $C_1$ field.
Using \cite{JCT}, the conjecture has been understood in the case when $\mr{char}(K) = \mr{char}(k)$.
However, little is known in the case of mixed characteristic i.e. $\mr{char}(K) \neq \mr{char}(k)$. 
In this note we prove the conjecture for the moduli space of stable locally free sheaves of fixed rank and determinant on a smooth, projective, geometrically connected curve defined over such a $C_1$ field.

\begin{note}\label{notin7}
Let $K$ be a field of characteristic $0$. 
Fix integers $r,d$ coprime with $r\geq 2$.
Let $X_K$ be a smooth, projective, geometrically connected curve of genus $g\geq 2$ defined over $K$. Assume there exists a line bundle $\mc{L}_K$ on $X_K$ of degree $d$. 
\end{note}

Let $M^s_{X_{K},\mc{L}_{K}}(r,d)$ be the moduli space of stable locally free sheaves of rank $r$ and determinant 
$\mc{L}_K$ (see Definition \ref{defimsfd} and Remark \ref{mklcorep}). Denote by  $M^s_{X_{\ov{K}},\mc{L}_{\ov{K}}}(r,d)$ 
the moduli space of stable locally free sheaves of  rank $r$ and determinant $\mc{L}_{\ov{K}}:= \mc{L}_K \otimes_{K} \ov{K}$ over the curve $X_{\ov{K}}:= X_K \times_K \msp(\ov{K})$. 
Since the functor $\mc{M}^{s}_{X_K,\mc{L}_{K}}(r,d)$ is universally corepresented by $M^{s}_{X_K,\mc{L}_{K}}(r,d)$, the moduli space $M^s_{X_{\ov{K}},\mc{L}_{\ov{K}}}(r,d)$ is isomorphic to $M^s_{X_{K},\mc{L}_{K}}(r,d) \times_{K} \msp(\ov{K})$.
By \cite{SFV}, $M^s_{X_{\ov{K}},\mc{L}_{\ov{K}}}(r,d)$ is a unirational variety and therefore rationally connected. 
Hence the moduli space $M^s_{X_{K},\mc{L}_{K}}(r,d)$ is a rationally connected variety.
Now suppose that $K$ is the fraction field of a Henselian discrete valuation ring with algebraically closed residue field.
The $C_1$ conjecture then predicts that $M^s_{X_{K},\mc{L}_{K}}(r,d)$ has a $K$-rational point.
In order to prove this, it suffices to show the existence of a stable locally free sheaf on $X_K$ of rank $r$ and determinant $\mc{L}_K$.

\vspace{0.2 cm}
The moduli of (semi)stable locally free sheaves of fixed rank and degree over a smooth, projective curve have been studied for decades and there is a plethora of results on the subject. However, for most of these results the curve is defined over an algebraically closed field. In fact when the field is not algebraically closed, there may not even exist invertible sheaves of certain degrees (see for example \cite{BI}, \cite{mes}). 
In this note we prove the following:

\begin{thm}[Theorem \ref{evb}, Corollary \ref{c1ex}]\label{thmin3}
Keep Notations \ref{notin7}.
There exists a stable locally free sheaf on $X_{K}$ of rank $r$ and 
determinant $\mc{L}_{K}$.

In particular, the moduli space of stable locally free sheaves over $X_{K}$ of rank $r$ and determinant $\mc{L}_K$ denoted $M_{X_K,\mc{L}_K}^{s}(r,d)$, has a $K$-rational point.   
\end{thm}

This result holds in much greater generality than $C_1$ fields and is therefore
of interest in its own right.
We use standard techniques from algebraic geometry to prove this theorem. Since the result holds for any field of characteristic $0$, it does not throw any light on the proof of the $C_1$ conjecture in the general case. Indeed it illustrates that even though at first sight the variety $M_{X_K,\mc{L}_K}^{s}(r,d)$ appears to be a good candidate for testing the conjecture in mixed characteristic, it is in fact a pathological example. 


\emph{Acknowledgements}: This paper answers the PhD question given to me by my supervisor Prof. H. Esnault.

\section{Main result}

We prove Theorem \ref{thmin3} stated in the introduction and show how it can be applied to the $C_1$ conjecture.

\vspace{0.2 cm}

\begin{thm}\label{evb}
Keep Notations \ref{notin7}.
There exists a geometrically stable locally free sheaf on $X_{K}$ of rank $r$ and 
determinant $\mc{L}_{K}$.
\end{thm}

\begin{proof}
  By \cite[Proposition $8.6.1$]{po} there exists a semistable, locally free sheaf of rank $r$ and 
 degree $d$ on $X_{\ov{K}}$. Since $\mr{Pic}^{0}(X_{\ov{K}})$ is an abelian variety and multiplication by $r$ is an isogeny, 
 one can show that there exists a semistable, locally free sheaf $\mc{E}_{\ov{K}}$ on $X_{\ov{K}}$
of rank $r$ and determinant $\mc{L}_{\ov{K}}$, where $\mc{L}_{\ov{K}}=\mc{L}_K \otimes \mo_{X_{\ov{K}}}$ is the base change of $\mc{L}_K$.
Furthermore, there exists an integer $b$ such that $\mc{E}_{\ov{K}} \otimes \mc{K}_{X_{\ov{K}}}^{\otimes b}$ is globally generated, 
where $\mc{K}_{X_{\ov{K}}}$ is the canonical divisor on $X_{\ov{K}}$. Since $X_{K}$ is a curve, this sheaf is a quotient of 
 $r+1$ copies of $\mo_{X_{\ov{K}}}$ i.e., we have the following surjective morphism:
 \[\bigoplus\limits_{i=1}^{r+1} \mo_{X_{\ov{K}}} \twoheadrightarrow \mc{E}_{\ov{K}} \otimes \mc{K}_{X_{\ov{K}}}^{\otimes b}.\]
 Since the determinant of $\mc{E}_{\ov{K}}$ is $\mc{L}_{\ov{K}}$, the kernel of this morphism is isomorphic to $\mc{L}_{\ov{K}}^\vee \otimes \mc{K}_{X_{\ov{K}}}^{-rb}$.
 In other words, $\mc{E}_{\ov{K}}$ is cokernel of a morphism \[\phi:\mc{L}_{\ov{K}}^\vee \otimes \mc{K}_{X_{\ov{K}}}^{-(r+1)b} \hookrightarrow \bigoplus\limits_{i=1}^{r+1}\mc{K}_{X_{\ov{K}}}^{-b}.\]
 Since $\mc{K}_{X_{\ov{K}}} \cong \mc{K}_{X_K} \otimes \mo_{X_{\ov{K}}}$, the sheaf $\mc{E}_{\ov{K}}$ is a $\ov{K}$-point of the affine space 
 \[U:=\mr{Hom}_{X_K}(\mc{L}_{{K}}^\vee \otimes \mc{K}_{X_{{K}}}^{-(r+1)b},\bigoplus\limits_{i=1}^{r+1}\mc{K}_{X_{{K}}}^{-b}).\]
 As local-freeness and semi-stability are open conditions, there exists a non-empty open subscheme $V$ of $U$ parameterizing those 
 homomorphisms whose cokernel is a semistable, locally free sheaf of rank $r$. Over any field of characteristic $0$, a nonempty Zariski open subset 
 of an affine space has a rational point.  Thus there exists a $K$-point of $V$ corresponding to a homomorphism
 such that the cokernel $\mc{F}_K$ is a locally free,  semistable rank $r$ sheaf with determinant $\mc{L}_K$ on $X_K$. 
 Since the degree of $\mc{L}_K$ is prime to $r$, $\mc{F}_K$ is also stable.
This proves the theorem.
 \end{proof}

Now we see an application of the above result.  

\vspace{0.2 cm}

\begin{defi}\label{defimsfd}
 Keep Notations \ref{notin7}.
Denote by $X_{T} := X_K \times_{\msp(K)} T$. 
We define a functor $\mc{M}_{X_K, \mc{L}_K}(r,d)$  as follows:
\[\mc{M}_{X_K, \mc{L}_K}(r,d) : \mr{Sch}^{\circ}/K \rightarrow \mr{Sets}\] 
\noindent such that for a $K$-scheme $T$,         
                                                                                      \[ \mc{M}_{X_K, \mc{L}_K}(r,d)(T):= \left\{ \begin{array}{l}
                                                                                  \mbox{ $S$-equivalence classes of } \mbox{ locally free sheaves } \mc{F} \mbox{ on  } X_{T}\\
                                                                                    \mbox{ such that for every geometric point } t \in T, \mc{F}_{t}  \mbox{ is a slope} \\ 
                                                                                   \mbox{ semistable sheaf of rank } r \mbox{ and   degree } d \mbox{ on } X_t \mbox{ and for}\\
                                                                                     \mbox{some invertible sheaf } \mc{Q} \mbox{ on } T, \mr{det}(\mc{F}) \simeq \pi^{*}_{X_{K}} \mc{L}_{K}\otimes \pi^{*}_{T}{\mc{Q}} \end{array} \right \} / \sim \]

\noindent where $\pi_{X_{K}}: X_{T} \rightarrow X_K$, $\pi_{T}: X_{T} \rightarrow T$ are the first and second projections respectively and 
$\mc{F} \sim \mc{F}'$ if and only if there exists an invertible sheaf $\mc{L}$ on $T$ such that $\mc{F} \simeq \mc{F}' \otimes \pi^{*}_{T}\mc{L}$.

\vspace{0.2 cm}
\noindent We denote by $\m^{s}_{X_K,\mc{L}_{K}}(r,d)$ the subfunctor for the stable sheaves.
Since $(r,d)$ are coprime and $X_K$ is integral, slope semistable sheaves are stable.
Hence $\m^{s}_{X_K,\mc{L}_{K}}(r,d)$ coincides with $\m_{X_K,\mc{L}_{K}}(r,d)$.
\end{defi}

\vspace{0.2 cm}

\begin{rem}\label{mklcorep}
Denote by $\m^{s}_{X_K}(r,d)$ the moduli functor of isomorphism classes of stable locally free sheaves of rank $r$ and degree $d$. By \cite[Theorem $0.2$]{LA1} this functor is universally corepresented by a projective $K$-scheme $M^s_{K}(r,d)$. Recall the Picard functor $\mc{P}ic_{X_K}$ and the natural transformation $\m^s_{X_K}(r,d) \rightarrow \mc{P}ic_{X_K}$ which is defined by taking the determinant of the locally free sheaves.
This induces the determinant morphism $\mr{det} : M^s_{X_K}(r,d) \rightarrow \mr{Pic}(X_{K})$, where $\mr{Pic}(X_{K})$ is the Picard group scheme of $X_K$.
Using the property of universal categorical quotients, one can show that $\mc{M}^{s}_{X_K,\mc{L}_{K}}(r,d)$ is universally corepresented by $\mr{det}^{-1}(\mc{L}_K)$ which we denote by $M^{s}_{X_K,\mc{L}_{K}}(r,d)$.
For a complete proof see \cite[Proposition $2.3$]{ink3} (replacing $R$ by $K$). Since $(r,d)$ are coprime in our current settings, we have $M^{s}_{X_K,\mc{L}_{K}}(r,d)$ is in fact a projective $K$-scheme.

 \end{rem}

\vspace{0.2 cm}

\begin{cor}\label{c1ex}
Keep Notations \ref{notin7}.
The moduli space $M_{X_K,\mc{L}_K}^{s}(r,d)$ has a $K$-rational point. 
\end{cor}

\begin{proof}
By Remark \ref{mklcorep}, $M_{X_K,\mc{L}_K}^{s}(r,d)$ corepresents the functor $\m^{s}_{X_K,\mc{L}_{K}}(r,d)$. 
By Theorem \ref{evb} there exists a stable locally free sheaf on $X_{K}$ of rank $r$ and determinant $\mc{L}_{K}$. 
Then by definition of corepresentability, 
the  moduli space $M_{X_K,\mc{L}_K}^{s}(r,d)$ has a $K$-rational point.
This proves the corollary.
\end{proof}

As mentioned in the introduction the variety $M^s_{X_K,\mc{L}_{K}}(r,d)$ is rationally connected.
The above corollary trivially proves the $C_1$ conjecture in mixed characteristic for this variety.

\vspace{0.2 cm}

\begin{cor}
Let $K$ be the fraction field of a Henselian discrete valuation ring with algebraically closed residue field $k$. 
Assume that the characteristic of $K$ is $0$ and that of $k$ is $p>0$. 
The $C_1$ conjecture holds for the variety $M^s_{X_K,\mc{L}_{K}}(r,d)$.
\end{cor}

\begin{proof}
By Corollary \ref{c1ex} the variety $M^s_{X_K,\mc{L}_{K}}(r,d)$ has a $K$-rational point.
\end{proof}


\begin{thebibliography}{1}

\bibitem{BI}
U.~Bhosle and I.~Biswas.
\newblock Stable real algebraic vector bundles over a {K}lein bottle.
\newblock {\em Transactions of the American Mathematical Society},
  360(9):4569--4595, 2008.

\bibitem{JCT}
J.~L. Colliot-Th{\'e}l{\`e}ne.
\newblock Vari{\'e}t{\'e}s presque rationnelles, leurs points rationnels et
  leurs d{\'e}g{\'e}n{\'e}rescences.
\newblock In {\em Arithmetic geometry}, pages 1--44. Springer, 2010.

\bibitem{ink3}
I.~Kaur.
\newblock Smoothness of moduli space of stable torsionfree sheaves with fixed
  determinant in mixed characteristic.
\newblock In {\em "Analytic and Algebraic Geometry"}, pages 173--186. Springer,
  2017.

\bibitem{L}
S.~Lang.
\newblock On quasi algebraic closure.
\newblock {\em Annals of Math.}, 55:373--390, 1952.

\bibitem{LA1}
A.~Langer.
\newblock Semistable sheaves in positive characteristic.
\newblock {\em Annals of Math}, 159:251--276, 2004.

\bibitem{po}
J.~Le~Potier.
\newblock {\em Lectures on vector bundles}, volume~54.
\newblock Cambridge University Press, 1997.

\bibitem{mes}
N.~Mestrano.
\newblock Conjecture de franchetta forte.
\newblock {\em Inventiones mathematicae}, 87(2):365--376, 1987.

\bibitem{SFV}
C.S. Seshadri.
\newblock {\em Fibres vectoriels sur les courbes algebriques}, volume 14023.
\newblock Ast{\'e}risque 96, Paris, 1982.

\end{thebibliography}
\end{document}